\documentclass[a4paper,12pt,wlayout]{article}

\usepackage{array}
\usepackage{amssymb}
\usepackage{amsmath}
\usepackage{amsbsy}
\usepackage{epsfig}
\usepackage{mathrsfs}
\usepackage{graphics}
\usepackage{graphicx}
\usepackage{pstricks}
\usepackage{float}
\usepackage{amsthm}
\usepackage{slashed}

\newcommand{\eps}{\varepsilon}
\newcommand{\lb}{\label}

\newcommand{\go}{\rightarrow}

\newcommand{\ee}{\end{equation}}
\newcommand{\be}{\begin{equation}}
\newcommand{\bea}{\begin{eqnarray}}
\newcommand{\eea}{\end{eqnarray}}
\newcommand{\sbea}{\begin{subequations}\begin{eqnarray}}
\newcommand{\seea}{\end{eqnarray}\end{subequations}}
\newcommand{\ees}{\end{equation*}}
\newcommand{\bes}{\begin{equation*}}
\newcommand{\beas}{\begin{eqnarray*}}
\newcommand{\eeas}{\end{eqnarray*}}

\newcommand{\rf}[1]{(\ref{#1})}

\newtheorem{theorem}{Theorem}[section]

\newenvironment{remark}[1][Remark]{\begin{trivlist}
		\item[\hskip \labelsep {\bfseries #1}]}{\end{trivlist}}

\begin{document}

\title{Asymptotic decay towards steady states of solutions to very fast and singular diffusion equations.}

\author{Georgy Kitavtsev\thanks{Middle East Technical University, Northern Cyprus Campus, Kalkanlı, Güzelyurt, KKTC, 10 Mersin, Turkey. {\it email: georgy.kitavtsev@gmail.com}} \and Roman M. Taranets\thanks{Institute of Applied Mathematics and Mechanics of the NASU, G.~Batyuka Str. 19, 84100, Sloviansk, Ukraine. {\it email: taranets\_r@yahoo.com}}}

\maketitle

\begin{abstract}
We analyze long-time behavior of solutions to a class of problems related to very fast and singular diffusion porous medium equations having nonhomogeneous in space and time source terms with zero mean. In dimensions two and three, we determine critical values of porous medium exponent for the asymptotic $H^1$-convergence of the solutions to a unique nonhomogeneous positive steady state generally to hold.
\end{abstract}

%%%%%%%%%%%%%%%%%%%%%%%%%%%%%%%%%%%%%%%%%%%%%%%%%%%%%%%%%%%%%%%%%
\section{Introduction and main results}
%%%%%%%%%%%%%%%%%%%%%%%%%%%%%%%%%%%%%%%%%%%%%%%%%%%%%%%%%%%%%%%%%

We study initial boundary value problem of the form:
\sbea
\lb{ME1}
v_t&=&v^k(\Delta v-f(x,t))\ \text{in}\ Q_T=\Omega\times(0,\,T),\\
\lb{ME2}
\nabla v\cdot\mathbf{n}&=&0\ \text{in}\ \partial\Omega\times(0,\,T),\\
\lb{ME3}
v(x,\,0)&=&v_0(x)\ge 0,
\seea
where $k>1$, $\Omega\subset\mathbb{R}^N$ is a bounded convex domain with Lipschitz boundary,
\be
\lb{f}
f\in L^s(0,+\infty;L^r(\Omega))\quad\text{and}\quad\int\limits_{\Omega} f(x,t)\,dx=0\quad\forall\, t\ge 0.
\ee
with some $s\in[1,\,+\infty]$ and $r\in(N/2,\,+\infty]$.
Note that \rf{ME1}--\rf{ME3} and \rf{f} together imply
\be
\lb{CM}
\displaystyle\int\limits_{\Omega}\tfrac{dx}{v^{k-1}}=\int\limits_{\Omega}\tfrac{dx}{v_0^{k-1}}=M\quad\text{for all}\ t>0.
\ee

The motivation to study \rf{ME1}--\rf{ME3} with source term \rf{f} comes from the following relation of this problem to the so called very fast and singular diffusion equations of the porous medium type (PME) (see~\cite{DP97}--\cite{GM13} and references therein). By denoting
\be
\lb{u_v}
u=\tfrac{1}{v^{k-1}}\quad\text{and}\quad m=\tfrac{k}{k-1}
\ee
system \rf{ME1}--\rf{ME3} transforms into the PME type initial boundary value problem
\sbea
\lb{PME1}
u_t&=&\mathrm{div}\bigl(\tfrac{\nabla u}{u^m}\bigr)+ \tfrac{1}{m-1} f(x,t)\ \text{in}\ Q_T=\Omega\times(0,\,T),\\
\lb{PME2}
\nabla u\cdot\mathbf{n}&=&0\ \text{in}\ \partial\Omega\times(0,\,T),\\
\lb{PME3}
u(x,\,0)&=&\tfrac{1}{v^{k-1}_0(x)}.
\seea
The latter problem has the source function $\frac{1}{m-1} f(x,t)$ and obeys the conservation law
\be
\lb{PCM}
\displaystyle\int\limits_{\Omega}u\,dx=M\quad\text{for all}\ t\ge 0.
\ee
Exploring the relation between problems \rf{ME1}--\rf{ME3} and \rf{PME1}--\rf{PME3} one finds that exponent
$k=2$ in the former model corresponds to $m=2$ in the latter with \rf{PME1}, which represents singular diffusion PME. Accordingly, for $k\in(2,+\infty)$ one has $m\in(1,\,2)$  with \rf{PME1}, which represent very fast diffusion PMEs. The limit $k\go+\infty$ formally corresponds to exponent $\frac{1}{u}$ in \rf{PME1} with $m=1$.  On the other hand, for $k\in(1,2)$ one has $m\in(2,\,+\infty)$ with  exponent of  \rf{PME1} being more degenerate then in the singular diffusion PME.

In dimension $N=1$ regularity and asymptotic convergence to a unique positive steady state of solutions to \rf{PME1}--\rf{PME3} was shown in~\cite{FKT18,KT20}. Up to our knowledge asymptotic properties and regularity of solutions to problem  \rf{PME1}--\rf{PME3} equipped with a source term \rf{f} in dimensions higher then one in the given range $m\in(1,\,+\infty)$ were only investigated for the case $m=2$ under additional condition $f(\cdot,t)\go 0$ as $t\go\infty$ (see~\cite{SSZ16} and references there in). In this study, we extend the results of~\cite{SSZ16,FKT18,KT20} in two ways. Firstly, we present new results about regularity and asymptotical decay properties of solutions to problem  \rf{PME1}--\rf{PME3} with general exponent $m\in(1,\,\infty)$ in dimensions $N=2$ and $N=3$. Secondly, we will considered general type source terms $f(x,t)$ as in \rf{f} that, in particular, include non-zero stationary $f(x,\,t)=f(x)$.

We also consider the unique positive solution $v_\infty\in W^{2,\,r}(\Omega)$ of the stationary problem for \rf{ME1}--\rf{ME3}:
\sbea
\lb{SME1}
\Delta v_\infty&=&f_\infty(x)\ \text{in}\ \Omega,\\
\lb{SME2}
\nabla v_\infty\cdot\mathbf{n}&=&0\ \text{in}\ \partial\Omega,\\
\lb{SME3}
\int\limits_{\Omega}\tfrac{dx}{v_\infty^{k-1}}&=&M,
\seea
where we assume that there exists a limiting source function $f_\infty(x)$, such that $f(\cdot,\,t)\go f_\infty(.)$ as $t\go\infty$ as stated in Theorems 1.1-1.2. Existence and uniqueness of such $v_\infty>0$ for $f_\infty\in L^r(\Omega)$
with $r>N/2$ is shown in Appendix A.2.

The main results of our article are stated in the following two theorems.
\begin{theorem}[time-homogeneous source]
Let $f(x,t)=f_\infty(x)=f(x)$ in \rf{f} with $r\in(N/2,\,+\infty]$.
Assume that $v_0(x)\ge 0$ and in the case $N=2$,
\be
\lb{ic_as2}
k\ge\tfrac{2r-1}{r-1}\quad\text{and}\quad\int\limits_{\Omega}\tfrac{1}{v_0^{k+\eps}}<\infty\quad\text{for an enough small}\quad\eps>0,
\ee
while in the case $N=3$,
\be
\lb{ic_as3}
k>\tfrac{5r-3}{2r-3}\quad\text{and}\int\limits_{\Omega}\tfrac{1}{v_0^\frac{3k}{2}}<\infty.
\ee
Then for any non-negative solution $v(x,\,t)$ to problem \rf{ME1}--\rf{ME3},
satisfying constraint \rf{CM} for all $t>0$, the following asymptotic decay to the unique positive solution  $v_\infty$ to the
stationary problem \rf{SME1}--\rf{SME3} holds:
\be
\lb{H1_decay}
\|v(\cdot,t)-v_\infty\|_{H^1(\Omega)}\le C_P\|\nabla(v_0-v_\infty)\|_2\exp(-\bar{C}t)\quad\text{for all }\,t>0,
\ee
in particular, $v(\cdot,t)\go v_\infty$ in $H^1(\Omega)$ as $t\go+\infty$. Positive constants $C_P$ and $C$ in \rf{H1_decay} depend on $k,\,r,\, |\Omega|,\,M,\,\|f\|_r$ and $v_0$ only.
\end{theorem}
\begin{theorem}[time-inhomogeneous source]
Let $r\in(N/2,\,+\infty]$ in \rf{f} and
\be
\lb{f_tau}
\int_0^\infty\|f_\tau\|_2\,d\tau<\infty
\ee
hold. Let
\be
\lb{f_conv}
\|f(\cdot,\,t)-f_\infty\|_r\go 0 \text{ as }  t\go+\infty,\quad\|f(\cdot,\,t)-f_0(.)\|_r\go 0 \text{ as } t\go 0.
\ee
In case $N=2$, let $s=\infty$ for $k\in[\tfrac{2r-1}{r-1},\,+\infty)$ and any $s\in[1,\,\frac{r}{r-(k-1)(r-1)})$
for $k\in(\frac{r}{r-1},\,\tfrac{2r-1}{r-1})$ in \rf{f}.

In case $N=3$, let $s=\infty$ for $k\in(\tfrac{5r-3}{2r-3},\,+\infty)$ and any $s\in[1,\tfrac{r(k+2)}{r(5-2k)+3(k-1)}]$ for $k\in[\frac{r}{r-3/2},\,\tfrac{5r-3}{2r-3}]$ in \rf{f}.

Additionally, assume that $v_0(x)\ge 0$, $E(v_0)<0$ in \rf{en_def}, and \rf{ic_as2} or \rf{ic_as3}
hold for $N=2$ or $N=3$, respectively.

Then any non-negative solution $v(x,\,t)$ to problem \rf{ME1}--\rf{ME3},
satisfying constraint \rf{CM} for all $t>0$, converges to the unique positive solution  $v_\infty$ to the
stationary problem \rf{SME1}--\rf{SME3} in $H^1(\Omega)$, i.e.
\be
\lb{non_stat_conv}
\|v(\cdot,t)-v_\infty\|_{H^1(\Omega)}\go 0 \text{ as }  t\go+\infty.
\ee
\end{theorem}
Theorems 1.1--1.2 indicate  for $N=2$ and $N=3$ critical values $k_{cr}(r)$ of exponent $k$ in \rf{ME1}--\rf{ME3}  (see also \rf{k_cr} and \rf{k_cr1} below) for the asymptotic $H^1$-convergence of dynamical solutions to a unique positive steady state generally to hold. Additionally, using transformation \rf{u_v} statements of Theorems 1.1--1.2 can be rewritten in terms of solutions to PME problems  \rf{PME1}--\rf{PME3}. In particular, \rf{H1_decay} and \rf{non_stat_conv} imply asymptotic decays
\bes
\|u^\frac{1}{1-k}(\cdot,t)-u_\infty^\frac{1}{1-k}\|_{H^1(\Omega)}\le C_P\|\nabla\bigl(u_0^\frac{1}{1-k}-u_\infty^\frac{1}{1-k}\bigr)\|_2\exp(-\bar{C}t)\quad\text{for all }\,t>0
\ees 	 
and
\bes
\|u(\cdot,t)^\frac{1}{1-k}-u_\infty^\frac{1}{1-k}\|_{H^1(\Omega)}\go 0 \text{ as }  t\go+\infty.
\ees
respectively, where $u_\infty$ is the solution of stationary problem to  \rf{PME1}--\rf{PME3}:
\beas
\mathrm{div}\bigl(\tfrac{\nabla u}{u^m}\bigr)&=&-\tfrac{1}{m-1}f_\infty(x)\ \text{in}\ \Omega,\\
\nabla u_\infty\cdot\mathbf{n}&=&0\ \text{in}\ \partial\Omega,\\
\int\limits_{\Omega}u_\infty\,dx&=&M,
\eeas
	
Let us also comment on two important applications of Theorems 1.1--1.2. For $k=2$ problem \rf{ME1}--\rf{ME3} corresponds to the singular heat equation \rf{PME1} with $m=2$ and in case $N=2$ further related by that to 3D viscous sheet model analyzed in~\cite{KT20}--\cite{KLN11}, the one that evoked primary motivation for this study. The source term considered in~\cite{FKT18,KT20} is homogeneous in time $f(x)\in L^\infty(\Omega)$. Theorem 1.1 for $N=2,\,k=2$ and $r=\infty$ combined with positivity and smoothness of dynamical solutions for all $t>0$ provided $v_0(x)>0$ (see Appendix A.1), then together imply asymptotic exponential decay formula \rf{H1_decay}. This result shows an exponential decay of the solutions to  the 3D viscous sheet model to a selected steady state, which depends only on the initial height and velocity distributions in the sheet. We also note that for $N=2$ and $r=+\infty$ the critical value of the exponent
\rf{k_cr} is $k_{cr}(+\infty)=2$. For $r=+\infty$ and $k\in(1,\,2)$ a convergence to a nonhomogeneous positive steady state then generally might not be expected, but Theorem 1.2 still identifies $s_{cr}=\tfrac{r(k+2)}{r(5-2k)+3(k-1)}=\tfrac{k+2}{5-2k}>1$ such that decay \rf{non_stat_conv} to constant state (see \rf{const_ss} in Appendix A.2) holds for any time exponent $s\in[1,\,s_{cr})$.

The second application comes from models closely related to 3D Penrose-Fife phase-field model~\cite{PF90} corresponding to $N=3$ and $k=2$ in problem \rf{ME1}--\rf{ME3} (see~\cite{SSZ16} and references therein). We note that the minimal regularity for this model suggested in~\cite{SSZ16}, i.e. $s=2$ and $r=3$ is exactly reflected in conditions of Theorem 1.2. Interestingly, here we land on another border case $k=\frac{r}{r-3/2}=2$ and $s=s_{cr}=\tfrac{r(k+2)}{r(5-2k)+3(k-1)}=2$. Additionally, for $N=3$ and $r=+\infty$ one has the critical value of the exponent \rf{k_cr1} $k_{cr}(+\infty)=2.5$ and Theorem 1.2 still identifies $s_{cr}(k)>1$ for any $k\in(1,\,2.5]$.

The organization of the article is as follows. Proofs of Theorems 1.1--1.2 are provided in section~3 for $N=2$ and section~4 for $N=3$. In section~2 we collect some basic estimates used often in the proofs. In the discussion section we state interesting open questions and outlook.

%%%%%%%%%%%%%%%%%%%%%%%%%%%%%%%%%%%%%%%%%%%%%%%%%%%%%%%%%%%%%%%%%
\section{Basic estimates}
%%%%%%%%%%%%%%%%%%%%%%%%%%%%%%%%%%%%%%%%%%%%%%%%%%%%%%%%%%%%%%%%%

By multiplying \rf{ME1} with $\Delta v-f(x,t)$ and integrating by parts one shows the energy equality
\be
\lb{EE}
\displaystyle\tfrac{d}{dt}E(v(t))+\int\limits_{\Omega} v^k(\Delta v-f)^2\,dx=\int\limits_{\Omega} vf_t\,dx
\ee
with
\be
\lb{en_def}
E(v)=\int\limits_{\Omega}[\tfrac{1}{2}|\nabla v|^2+fv]\,dx,
\ee
that consequently implies
\be
\lb{grad_v}
\|\nabla v(t)\|_2\le C \text{ for all }  t>0
\ee
provided (\ref{f_tau}) holds, with constant $C$ depending on the value of integral in \rf{f_tau} as well as on b$E(v_0),\,\|f_0\|_2,\,M, \,|\Omega|$. Proof of \rf{grad_v} is given in Appendix~B.

Let us denote the difference between dynamical and stationary solutions by $w=v(x,t)-v_\infty(x)$. Then one
can rewrite \rf{ME1} as follows
\bes
v_t=v^k(\Delta v-f(x,t))=v^k\Delta w+v^k(f_\infty(x)-f(x,t)).
\ees
Multiplying the last equation by $- \Delta w$ and integrating by parts results in
\be
\lb{w_es}
\tfrac{1}{2}\tfrac{d}{dt}\int\limits_{\Omega}|\nabla w|^2\,dx+\int\limits_{\Omega} v^k|\Delta w|^2\,dx=
\int\limits_{\Omega} v^k(f(x,t)-f_\infty(x))\Delta w\,dx.
\ee
Additionally, we will use the following estimate
\be
\lb{diss_est}
\int\limits_{\Omega}|\nabla w|^2\,dx\le C_{S,l}^2\Bigl(\int\limits_{\Omega} v^k|\Delta w|^2\,dx\Bigr)\Bigl(\int\limits_{\Omega} v^{\frac{kl}{2-l}}\,dx\Bigr)^{\frac{l-2}{l}},
\ee
which holds for any $l>2$ in the case $N=2$ and for $l\in(2,\,6]$ in the case $N=3$ with $C_{S,l}$ being the constant in Poincare-Sobolev embedding
\bes
\|w-\bar{w}\|_l\le C_{S,l}\|\nabla w\|_2,\quad\text{where}\quad\bar{w}=\tfrac{1}{|\Omega|}\int\limits_{\Omega}w\,dx.
\ees
Combining the latter estimate with the following calculation, taking into account $\nabla w\cdot\mathbf{n}=0$ on $\partial\Omega$,
\beas
\int\limits_{\Omega}|\nabla w|^2\,dx&=&-\int\limits_{\Omega} (w - \bar{w})\Delta w\,dx \le
\Bigl ( \int\limits_{\Omega} v^{\frac{k}{2}}|\Delta w|^2\,dx \Bigr)^{\frac{1}{2}}
\|w - \bar{w} \|_l \|v^{-\frac{k}{2}}\|_{\frac{2l}{l-2}}\nonumber\\
&\le&C_{S,l} \Bigl ( \int\limits_{\Omega} v^{\frac{k}{2}}|\Delta w|^2\,dx \Bigr)^{\frac{1}{2}} \|\nabla w\|_2
\Bigl(\int\limits_{\Omega} v^{\frac{kl}{2-l}}\,dx \Bigr)^{\frac{l-2}{2l}},
\eeas
which implies \rf{diss_est}.

Finally, multiplying \rf{ME1} by $v^{-p-1}$ and integrating by parts one obtains
\be
\lb{Entrops}
\tfrac{d}{dt}\int\limits_{\Omega}\tfrac{dx}{v^p}+\tfrac{4(p+1-k)p}{(p-k)^2}
\int\limits_{\Omega} \bigl|\nabla\bigl(\tfrac{1}{v^{(p-k)/2}}\bigr)\bigr|^2\,dx= p\int\limits_{\Omega}\tfrac{f\,dx}{v^{p+1-k}}
\ee
holding for any $p>k-1>0$ and $p \neq k$.

%%%%%%%%%%%%%%%%%%%%%%%%%%%%%%%%%%%%%%%%%%%%%%%%%%%%%%%%%%%%%%%%%
\section{Two-dimensional case}
%%%%%%%%%%%%%%%%%%%%%%%%%%%%%%%%%%%%%%%%%%%%%%%%%%%%%%%%%%%%%%%%%
%%%%%%%%%%%%%%%%%%%%%%%%%%%%%%%%%%%%%%%%%%%%%%%%%%%%%%%%%%%%%%%%%
\subsection{Time-homogeneous case}
%%%%%%%%%%%%%%%%%%%%%%%%%%%%%%%%%%%%%%%%%%%%%%%%%%%%%%%%%%%%%%%%%
In this section, we collect and prove results on exponential asymptotic decay to the steady state $v_\infty(x)$ of solutions
to problem \rf{ME1}--\rf{ME3} considered with source function \rf{f} being independent of time,
i.e. $f(x,t)=f_\infty(x)=f(x)$. It turns out that given $f\in L^r(\Omega)$ with $r\in(1,\,+\infty]$ there exists a critical value of exponent
\be
\lb{k_cr}
k_{cr}(r)=\tfrac{2r-1}{r-1}\ge 2,
\ee
such that for $k\ge k_{cr}(r)$ we can show an exponential asymptotic decay estimate \rf{H1_decay}
for any given positive $|\Omega|,M,\|f\|_r$ and initial condition $v_0(x) \ge0$, provided \rf{ic_as2} holds.

Showing unconditional asymptotic decay in case $k<k_{cr}(r)$ remains then an open problem. Note, that given $f\in L^\infty(\Omega)$ one has $r=+\infty$ and $k_{cr}(+\infty)=2$. 

Below, in subsections 3.1.1 and 3.1.2, we derive uniform in time bounds on $\|\frac{1}{v}\|_\frac{kl}{l-2}$ for certain $l>2$.
Combining them with \rf{diss_est} and $f_\infty(x)=f(x)$ then will yield the estimate
\be
\lb{Msc_es}
\tfrac{1}{2}\tfrac{d}{dt}\int\limits_{\Omega}|\nabla w|^2\,dx+\bar{C}\int\limits_{\Omega}|\nabla w|^2\,dx\le0\quad\text{for}\quad t>0
\ee
and some constant $\bar{C}>0$.
Applying comparison principle to the last in\-equa\-li\-ty will show
\bes
\int\limits_{\Omega}|\nabla w|^2\,dx\le\int\limits_{\Omega}|\nabla w_0|^2\,dx\exp(-2\bar{C}t),
\ees
whence substituting $w=v(x,t)-v_\infty(x)$ gives
\bes
\|\nabla(v(\cdot,t)-v_\infty)(.)\|_2\le\|\nabla(v_0-v_\infty)\|_2\exp(-\bar{C}t).
\ees
The last estimate combined together with Poincare inequality, conversations \rf{CM} and \rf{SME3}, non-negativity of $v(x,t)$ and positivity of  $v_\infty(x)\in C^\lambda(\bar{\Omega})$   (see Appendix A.2) will, finally, imply \rf{H1_decay} with $C_P$ being Poincare's constant.

%%%%%%%%%%%%%%%%%%%%%%%%%%%%%%%%%%%%%%%%%%%%%%%%%%%%%%%%%%%%%%%%%
\subsubsection{Case $k>k_{cr}(r)$}
%%%%%%%%%%%%%%%%%%%%%%%%%%%%%%%%%%%%%%%%%%%%%%%%%%%%%%%%%%%%%%%%%
First, let us estimate the last term in \rf{Entrops} using H\"{o}lder inequality as
\be
\lb{RHS}
\int\limits_{\Omega}\tfrac{f\,dx}{v^{p+1-k}}\le\|f\|_r\Bigl(\int\limits_{\Omega}\tfrac{dx}{v^\frac{(p+1-k)r}{r-1}}\Bigr)^{\frac{r-1}{r}}.
\ee
Conservation \rf{CM} implies existence of $x_0(t)\in\Omega$, such that $v(t,\,x_0(t))=(|\Omega|/M)^{1/(k-1)}$. Combining this fact with Sobolev-Poincare inequality one estimates
\be
\lb{MTest}
\int|\nabla\Bigl(\tfrac{1}{v^{(p-k)/2}}\Bigr)|^2\,dx\ge\tfrac{1}{C_{S,q}^2}\Bigl(\int\limits_{\Omega}\Bigl|\tfrac{1}{v^{(p-k)/2}}-\Bigl(\tfrac{M}{|\Omega|}\Bigr)^\frac{p-k}{2k-2}\Bigr|^q\,dx\Bigr)^{2/q}
\ee
for any $q\in(1,\,\infty)$, where $C_{S,q}$ denotes Sobolev constant that depends on $|\Omega|$ and choice of $q$.
Using inequalities $\|v-w\|_q\ge|\|v\|_q-\|w\|_q|$ and
\bes
(a-b)^2\ge\tfrac{a^2}{2}-b^2,
\ees
estimates \rf{Entrops}, \rf{RHS}--\rf{MTest} together imply
\be
\lb{Entrops1}
\tfrac{d}{dt}\int\limits_{\Omega}\tfrac{dx}{v^p}+\tfrac{2(p+1-k)p}{C_{S,q}^2(p-k)^2}\Bigl(\int\limits_{\Omega}\tfrac{1}{v^\frac{q(p-k)}{2}}\,dx\Bigr)^\frac{2}{q}
\le\tfrac{4(p+1-k)p}{C_{S,q}^2(p-k)^2}\Bigl(\tfrac{M}{|\Omega|}\Bigr)^{\frac{p-k}{k-1}}|\Omega|^{\frac{2}{q}}
+p\|f\|_r\Bigl(\int\limits_{\Omega}\tfrac{dx}{v^\frac{(p+1-k)r}{r-1}}\Bigr)^{\frac{r-1}{r}}_.
\ee
Let us set
\be
\lb{pdef}
p=\tfrac{(k-1)(2r-1)}{r},\quad q=\tfrac{2p}{p-k}>1,
\ee
and denote
\bes
y(t)=\int\limits_{\Omega}\tfrac{dx}{v^p}.
\ees
Due to definition \rf{pdef} and constraint $k>k_{cr}(r)$, one has $p>k$ and \rf{Entrops1} becomes
\bea
\lb{Entrops2}
&&\hspace{-0.cm}y'(t)+\tfrac{2(k-1)^2(2r-1)(r-1)}{C_{S,\frac{2(k-1)(2r-1)}{(k-2)r-k+1}}^2((k-2)r-k+1)^2}y(t)^\frac{(k-2)r-k+1}{(k-1)(2r-1)}\le\tfrac{(k-1)(2r-1)}{r}\|f\|_rM^{\frac{r-1}{r}}\nonumber\\
&&\hspace{-0.cm}+\tfrac{4(k-1)^2(2r-1)(r-1)}{C_{S,\frac{2(k-1)(2r-1)}{(k-2)r-k+1}}^2((k-2)r-k+1)^2}
\bigl(\tfrac{M}{|\Omega|}\bigr)^\frac{(k-2)r-k+1}{(k-1)r}|\Omega|^{\frac{(k-2)r-k+1}{(k-1)(2r-1)}},
\eea
where in the first term at the right-hand side we used \rf{CM}. Hence, the right hand side turns out to be a constant depending on given values $k,r,|\Omega|,M,\|f\|_r$.  Consequently, the last inequality implies a uniform in time bound
\be
\lb{p_bound}
y(t)=\int\limits_{\Omega}\tfrac{dx}{v^\frac{(k-1)(2r-1)}{r}}\le C=\max\{y(0),\,C_0\}
\ee
with constant
\bes
C_0=\Bigl[2\bigl(\tfrac{M}{|\Omega|}\bigr)^\frac{(k-2)r-k+1}{(k-1)r}|\Omega|^\frac{(k-2)r-k+1}{(k-1)(2r-1)}+\|f\|_rM^\frac{r-1}{r}\tfrac{C_{S,\frac{2(k-1)(2r-1)}{(k-2)r-k+1}}^2((k-2)r-k+1)^2}{2(k-1)r(r-1)}\Bigr]^\frac{(k-1)(2r-1)}{(k-2)r-k+1}_.
\ees
Next, choosing
\be
\lb{l}
l=\tfrac{2(k-1)(2r-1)}{(k-2)r-k+1}\in(2,\infty),\quad\text{i.e. so that}\quad\tfrac{(k-1)(2r-1)}{r}=\tfrac{kl}{l-2},
\ee
in estimates \rf{w_es}--\rf{diss_est} and using \rf{p_bound} yields \rf{Msc_es} with
\bes
\bar{C}=C_{S,\frac{2(k-1)(2r-1)}{(k-2)r-k+1}}^{-2}C^{-\frac{kr}{(k-1)(2r-1)}},
\ees
and, consequently, decay estimate \rf{H1_decay} holds.

We note that the derived exponential asymptotic rate $\bar{C}$ degenerates in the limit $k\go k^+_{cr}(r)$, namely,
\be
\lb{as_deg}
C_{S,\frac{2(k-1)(2r-1)}{(k-2)r-k+1}}^{-2}C^{-\frac{kr}{(k-1)(2r-1)}}\go 0\quad\text{as}\quad k\go k^+_{cr}(r).
\ee
To justify this we observe that
\bes
(k-2)r-k+1\go 0^+ \text{ and }  C^{\frac{kr}{(k-1)(2r-1)}}\sim2^\frac{kr}{(k-2)r-k+1}\go +\infty \text{ as }  k\go k^+_{cr}(r).
\ees
and by  Theorem 3.2 of~\cite{MTSO17} (see also Theorem 8.5(ii) of~\cite{LL01})
\bes
C_{S,\frac{2(k-1)(2r-1)}{(k-2)r-k+1}}\sim\tfrac{C(|\Omega|,r)}{\sqrt{(k-2)r-k+1}}\go +\infty  \text{ as }  k\go k^+_{cr}(r).
\ees

%%%%%%%%%%%%%%%%%%%%%%%%%%%%%%%%%%%%%%%%%%%%%%%%%%%%%%%%%%%%%%%%%
\subsubsection{Case $k=k_{cr}(r)$}
%%%%%%%%%%%%%%%%%%%%%%%%%%%%%%%%%%%%%%%%%%%%%%%%%%%%%%%%%%%%%%%%%
For any parameter $\beta\in(0,1)$ one can estimate the last term in \rf{Entrops} using H\"{o}lder inequality and \rf{CM} as
\be
\lb{RHS1}
\int\limits_{\Omega}\tfrac{f\,dx}{v^{p+1-k}}=\int\limits_{\Omega}\tfrac{f\,dx}{v^{p+1-k-\beta}\cdot v^\beta}\le\|f\|_rM^{\frac{\beta}{k-1}}\Bigl(\int\limits_{\Omega}\tfrac{dx}{v^{(p+1-k-\beta)k'}}\Bigr)^{\frac{1}{k'}},
\ee
where
\bes
k'=\tfrac{(k-1)r}{(r-1)(k-1)-r\beta}\ \underbrace{=}_{k=k_{cr}(r)}\ \tfrac{k-1}{1-\beta}>1.
\ees
Let us choose $\beta$ so that $(p+1-k-\beta)k'=p$, i.e.
\be
\lb{beta}
\beta=\tfrac{(k-1)^2-p(k-2)}{p-k+1}.
\ee
Condition $\beta\in(0,1)$ implies the following constraint on $p$
\be
\lb{p_constr}
p\in \bigl(k,\tfrac{(k-1)^2}{k-2}\bigr)
\ee
with the latter interval having non-zero length $1/(k-2)$.

Next, using \rf{RHS1}, \rf{beta}, again Sobolev-Poincare embedding \rf{MTest} with
\bes
q=\tfrac{2p}{p-k}>1
\ees
and by denoting
\bes
y(t)=\int\limits_{\Omega}\tfrac{dx}{v^p},
\ees
one obtains an analogue of inequality \rf{Entrops1}:
\be
\lb{Entrops3}
y'(t)+\tfrac{2(p+1-k)p}{C_{S,\frac{2p}{p-k}}^2(p-k)^2}y(t)^\frac{p-k}{p}\le\tfrac{4(p+1-k)p}{C_{S,\frac{2p}{p-k}}^2(p-k)^2}
\bigl(\tfrac{M}{|\Omega|}\bigr)^{\frac{p-k}{k-1}}|\Omega|^{\frac{p-k}{p}}
+p\|f\|_rM^\frac{1+(p-k)(2-k)}{(p-k+1)(k-1)}y(t)^{\frac{p-k}{p-k+1}}.
\ee
We rewrite \rf{Entrops3} in the following form
\be
\lb{Entrops4}
y'(t)\le C_1 z^\frac{p}{p-k+1}(t) - C_2 z(t)+C_3,
\ee
where we introduced the following notation:
\beas
z(t)&=&y(t)^\frac{p-k}{p},\  C_1=p\|f\|_rM^\frac{1+(p-k)(2-k)}{(p-k+1)(k-1)},\\[1.5ex] C_2&=&\tfrac{2(p+1-k)p}{C_{S,\frac{2p}{p-k}}^2(p-k)^2},\quad
C_3=2C_2\bigl(\tfrac{M}{|\Omega|}\bigr)^{\frac{p-k}{k-1}}|\Omega|^{\frac{p-k}{p}}.
\eeas
One observes that constants $C_1,C_2,C_3$ are positive if $p$ belongs to the range \rf{p_constr}.
The polynomial at the right-hand side of \rf{Entrops4} has two roots $\lambda_\pm$:
\be
\lb{lambda}
0<\lambda_- <\bigl(\tfrac{C_2(p-k+1)}{pC_1}\bigr)^\frac{p-k+1}{k-1}<\lambda_+
\ee
provided
\bes
C_3 < \tfrac{k-1}{p}\bigl(\tfrac{p-k+1}{pC_1}\bigr)^\frac{p-k+1}{k-1}C_2^\frac{p}{k-1}.
\ees
The latter condition can be rewritten as
\be
\lb{cond}
\tfrac{2p}{k-1}\bigl(\tfrac{M}{|\Omega|}\bigr)^{\frac{p-k}{k-1}}|\Omega|^{\frac{p-k}{p}}\bigl(\tfrac{pC_1}{p-k+1}\bigr)^\frac{p-k+1}{k-1}<C_2^\frac{p-k+1}{k-1}.
\ee
Observe that
\bes
\mathop {\lim} \limits_{p\go k^+} \tfrac{2p}{k-1}\bigl(\tfrac{M}{|\Omega|}\bigr)^{\frac{p-k}{k-1}}|\Omega|^{\frac{p-k}{p}}\bigl(\tfrac{pC_1}{p-k+1}\bigr)^\frac{p-k+1}{k-1}
=\tfrac{2k}{k-1}(k^2\|f\|_rM^\frac{1}{k-1})^\frac{1}{k-1} < +\infty,
\ees
while
\bes
\mathop {\lim} \limits_{p\go k^+}  C_2^\frac{p-k+1}{k-1} = +\infty,
\ees
where we use asymptotics
\bes
C_{S,\frac{2p}{p-k}}\sim\tfrac{C(|\Omega|,k)}{\sqrt{p-k}}\go\infty \quad\text{as}\quad p\go k^+,
\ees
see Theorem 3.2 of~\cite{MTSO17} (also Theorem 8.5(ii) of~\cite{LL01}).

Therefore, for $p$ being sufficiently close to value $k$ from above (depending on the given values $M,\,|\Omega|,\,\|f\|_r$ and $k=k_{cr}(r)$) inequality \rf{cond} and, consequently, \rf{lambda} hold with $\lambda_+\go+\infty$ as $p\go k^+$.

Now, to obtain a uniform bound on $y(t)$ for all $t>0$ we choose $p$ being sufficiently close to value $k$ from above, such that \rf{cond} holds as well as
\be
y(0)\le\lambda_+^\frac{p}{p-k}.
\ee
Such choice of $p$ is possible, because $y(0)$ stays bounded for $p$ being close to $k$ provided \rf{ic_as2} holds.
Then from \rf{Entrops4} and \rf{lambda} it follows that
\be
\lb{p_bound1}
y(t)=\int\limits_{\Omega}\tfrac{dx}{v^p}\le C=\min\Bigl\{y(0),\lambda_+^\frac{p}{p-k}\Bigr\}\quad\forall\, t>0.
\ee
Finally, fixing $p>k$ as above and choosing
\bes
l=\tfrac{2p}{p-k}\in(2,+\infty),\quad\text{i.e. so that}\quad p=\tfrac{kl}{l-2},
\ees
in estimates \rf{w_es}--\rf{diss_est} together with \rf{p_bound1} gives \rf{Msc_es} with
\bes
\bar{C}=C_{S,\frac{2p}{p-k}}^{-2}C^{-\frac{k}{p}},
\ees
and, consequently, decay estimate \rf{H1_decay} holds.

We note that the derived exponential asymptotic rate $\bar{C}$ degenerates in the limit $p\go k^+(r)$, namely
\bes
C_{S,\frac{2p}{p-k}}^{-2}C^{-\frac{k}{p}}\go 0\quad\text{as}\quad p\go k^+.
\ees
This follows from unboundedness of $\lambda_+$ in \rf{lambda} and $C_{S,\frac{2p}{p-k}}$ as $p\go k^+$.

%%%%%%%%%%%%%%%%%%%%%%%%%%%%%%%%%%%%%%%%%%%%%%%%%%%%%%%%%%%%%%%%%
\subsection{Time-inhomogeneous case}
%%%%%%%%%%%%%%%%%%%%%%%%%%%%%%%%%%%%%%%%%%%%%%%%%%%%%%%%%%%%%%%%%
In this section, we collect and prove results on asymptotic decay to the steady state $v_\infty(x)$ of solutions to problem \rf{ME1}--\rf{ME3} considered with time-dependent source function $f(x,\,t)$ in \rf{f}. We separate cases $k\ge k_{cr}(r)$ and $k\in(\frac{r}{r-1},\,k_{cr}(r))$ with $k_{cr}(r)$ being defined in \rf{k_cr}. It turns out that in the first case one can take $s=\infty$ in \rf{f} (similar to to the time-homogeneous case $f(x,t)=f(x)$), while in case $k\in(\frac{r}{r-1},\,k_{cr}(r))$ for a given space integrability $r\in(1,\infty]$ we can identify a maximal possible $s_{cr}=s_{cr}(k,r)\in[1,\infty)$, such that asymptotic decay to the steady state $v_\infty(x)$ is shown for any time integrability $s\in[1,\,s_{cr}(k,r))$ in \rf{f}. For $k\le\frac{r}{r-1}$ possibility for an asymptotic decay to  $v_\infty(x)$ remains then an open problem. 

We start with the following estimation of the right-hand side of \rf{w_es}
\bes
\int\limits_{\Omega} v^k(f-f_\infty)\Delta w\,dx\le\tfrac{1}{2}\int\limits_{\Omega} v^k(f-f_\infty)^2\,dx+\tfrac{1}{2}\int\limits_{\Omega}v^k|\Delta w|^2\,dx.
\ees
Combined with \rf{w_es} this implies
\be
\lb{w_es1}
\tfrac{d}{dt}\int\limits_{\Omega}|\nabla w|^2\,dx+\int\limits_{\Omega} v^k|\Delta w|^2\,dx\le
\int\limits_{\Omega} v^k(f-f_\infty)^2\,dx.
\ee
For any real $a>1$ one can estimate the right-hand side of \rf{w_es1} as
\bes
\int\limits_{\Omega} v^k(f-f_\infty)^2\,dx\le\|v\|_{ak}^k\|f-f_\infty\|_{\frac{2a}{a-1}}^2 .
\ees
Using Poincare-Sobolev embedding and \rf{CM} one shows
\bes
\|v\|_{ak}^k\le \bigl(C_{S,ak}\|\nabla v\|_2+\bigl(\tfrac{|\Omega|}{M}\bigr)^\frac{1}{k-1}|\Omega|^\frac{1}{ak}\bigr)^k.
\ees
Choosing $a=1+2/\eps$ with $\eps>0$ small enough one derives from the last two inequalities
\be
\lb{f_diff}
\int\limits_{\Omega}v^k(f-f_\infty)^2\,dx\le\tfrac{\bar{C}}{\eps^{k/2}}\|f-f_\infty\|_{2+\eps}^2,
\ee
where we use \rf{grad_v} and asymptotics
\bes
C_{S,1+2/\eps}\sim\tfrac{C(|\Omega|)}{\sqrt{\eps}} \text{ as }  \eps\go 0.
\ees
Finally, combining \rf{f_diff} with \rf{w_es1} gives
\be
\lb{w_es2}
\tfrac{d}{dt}\int\limits_{\Omega}|\nabla w|^2\,dx+\int\limits_{\Omega} v^k|\Delta w|^2\,dx\le
\tfrac{\bar{C}}{\eps^{k/2}}\|f-f_\infty\|_{2+\eps}^2.
\ee
Below, in subsections 3.2.1 and 3.2.2, we derive uniform in time bounds on $\|\frac{1}{v}\|_\frac{kl}{l-2}$ for certain $l>2$.
Combining them with \rf{diss_est} and \rf{w_es2} will yield asymptotic decay estimate \rf{H1_decay} for any $k>1$ and $r>2$. We note that this proof can be also generalized for the case $r\in(1,\,2]$ if one treats \rf{w_es} similarly as in section 4.2 for $N=3$.

%%%%%%%%%%%%%%%%%%%%%%%%%%%%%%%%%%%%%%%%%%%%%%%%%%%%%%%%%%%%%%%%%
\subsubsection{Case $k\ge k_{cr}(r)$}
%%%%%%%%%%%%%%%%%%%%%%%%%%%%%%%%%%%%%%%%%%%%%%%%%%%%%%%%%%%%%%%%%

Here, w.l.o.g., we assume that $s=\infty$ and $r\in(2,\,+\infty]$ in \rf{f}. In the first part of the proof, we follow exactly the lines of subsection 3.1.1 (in the case $k>k_{cr}(r)$) or subsection 3.1.2 (in the case $k=k_{cr}(r)$) establishing uniform in time bounds \rf{p_bound} or \rf{p_bound1}, respectively, i.e.
\bes
\int\limits_{\Omega}\tfrac{dx}{v^p}\le C\quad\forall t>0\quad\text{and fixed}\quad p\in(k,\,+\infty),\ C>0.
\ees
We note just that by deriving bounds \rf{p_bound} or \rf{p_bound1} $\|f\|_r$ should be now replaced by
$\|f\|_{L^\infty(0,+\infty;L^r(\Omega))}$.

Next, the latter bound combined with \rf{w_es2} and \rf{diss_est} taken with $l=\frac{2p}{p-k}$ implies
\bes
\tfrac{d}{dt}\int\limits_{\Omega}|\nabla w|^2\,dx+C_{S,\frac{2p}{p-k}}^{-2}C^{-\frac{k}{p}}\int\limits_{\Omega}|\nabla w|^2\,dx\le
\tfrac{\bar{C}}{\eps^{k/2}}\|f-f_\infty\|_{2+\eps}^2.
\ees
Applying comparison principle to the last inequality one deduces that
\bea
\int\limits_{\Omega}|\nabla w|^2\,dx&\le& \exp\Bigl\{-C_{S,\frac{2p}{p-k}}^{-2}C^{-\frac{k}{p}}\,t\Bigr\}\Big(\int\limits_{\Omega}|\nabla w_0|^2\,dx+\nonumber\\
&+&\tfrac{\bar{C}}{\eps^{k/2}}\int_0^t\exp\Bigl\{C_{S,\frac{2p}{p-k}}^{-2}C^{-\frac{k}{p}}\,\tau\Bigr\}\|f(\cdot,\tau)-f_\infty(.)\|_{2+\eps}^2\,d\tau\Big).
\eea
The last estimate would imply asymptotic decay
\be
\lb{as_decay2}
\int\limits_{\Omega}|\nabla w|^2\,dx\go 0\quad\text{as}\quad t\go+\infty
\ee
provided
\bes
\mathop {\lim} \limits_{t\go+\infty} \exp\Bigl\{-C_{S,\frac{2p}{p-k}}^{-2}C^{-\frac{k}{p}}\,t\Bigr\}\int_0^t\exp\Bigl\{C_{S,\frac{2p}{p-k}}^{-2}C^{-\frac{k}{p}}\,
\tau\Bigr\}\|f(\cdot,\tau)-f_\infty(.)\|_{2+\eps}^2\,d\tau=0,
\ees
which is true by \rf{f_conv}.
Using L'Hospital's rule one observes that the last limit holds if and only if
\bes
\lim_{t\go\infty}\|f(\cdot,t)-f_\infty(.)\|_{2+\eps}^2=0.
\ees
In that case \rf{as_decay2} holds and, after substituting $w=v(x,t)-v_\infty(x)$, one has
\bes
\|\nabla(v(\cdot,t)-v_\infty(.))\|_2\go 0 \text{ and } \|v(\cdot,t)-v_\infty(.)\|_{H^1(\Omega)}\go 0 \text{ as }  t\go+\infty.
\ees
We note that the decay rate here is determined by the slowest one between two: given convergence rate of $\|f(\cdot,\tau)-f_\infty(.)\|_{2+\eps}$ to zero and exponential $C_{S,\frac{2p}{p-k}}^{-2}C^{-\frac{k}{p}}$ one.

%%%%%%%%%%%%%%%%%%%%%%%%%%%%%%%%%%%%%%%%%%%%%%%%%%%%%%%%%%%%%%%%%
\subsubsection{Case $\tfrac{r}{r-1}<k< k_{cr}(r)$}
%%%%%%%%%%%%%%%%%%%%%%%%%%%%%%%%%%%%%%%%%%%%%%%%%%%%%%%%%%%%%%%%%

Here, we assume $f\in L^s(0,+\infty;L^r(\Omega))$ with $s\in[1,\infty)$ and $r\in(2,\,+\infty]$.
First, for any parameter $\beta\in[0,k-1)$ one estimates the last term in \rf{Entrops}, using H\"{o}lder inequality and \rf{CM}, as
\be
\lb{RHS2}
\int\limits_{\Omega}\tfrac{f\,dx}{v^{p+1-k}}=\int\limits_{\Omega}\tfrac{f\,dx}{v^{p+1-k-\beta}\cdot v^\beta}\le\|f\|_rM^{\tfrac{\beta}{k-1}}\Bigl(\int\limits_{\Omega}\tfrac{dx}{v^{(p+1-k-\beta)k'}}\Bigr)^{\frac{1}{k'}},
\ee
where
\be
\lb{k'}
k'=\tfrac{(k-1)r}{(r-1)(k-1)-r\beta}>1.
\ee
Let us choose $\beta$ so that $(p+1-k-\beta)k'=p$, i.e.
\be
\lb{beta1}
\beta=\tfrac{(k-1)^2r-p(k-1)}{r(p-k+1)}.
\ee
Condition $\beta\in[0,k-1)$ implies the following constraints on $p$:
\be
\lb{p_constr1}
p\in\bigl(\tfrac{2(k-1)r}{r+1},\,(k-1)r\bigr]
\ee
with the latter interval having non-zero length for $r>1$. Additionally, we demand $p+1-k-\beta>0$,
which, in view of \rf{beta1}, implies
\bes
p>\tfrac{(k-1)(2r-1)}{r}.
\ees
We note that for $1<k<k_{cr}(r)$ and $r>1$ one has
\bes
\tfrac{2(k-1)r}{r+1}<\tfrac{(k-1)(2r-1)}{r}<k,
\ees
and, therefore, we update the constraints on $p$ in \rf{p_constr1} to
\be
\lb{p_constr2}
p\in\bigl(k,\,(k-1)r\bigr]
\ee
with the latter interval having non-zero length, because $\frac{r}{r-1}<k$.
We conclude that for any $p$ taken from \rf{p_constr2} by \rf{Entrops} and \rf{RHS2}--\rf{beta1}
\bes
\tfrac{d}{dt}\int\limits_{\Omega}\tfrac{dx}{v^p}\le p\|f\|_rM^\frac{(k-1)r-p}{r(p-k+1)}\Bigl(\int\limits_{\Omega}\tfrac{dx}{v^p}\Bigr)^\frac{pr-(k-1)(2r-1)}{r(p-k+1)}
\ees
holds. Applying Gr\"{o}nwall inequality to this estimate results in the following bound
\be
\lb{p_bound2}
\int\limits_{\Omega}\tfrac{dx}{v^p}\le \Bigl[\Bigl(\int\limits_{\Omega}\tfrac{dx}{v_0^p}\Bigr)^\frac{(k-1)(r-1)}{r(p-k+1)}
+\tfrac{p(k-1)(r-1)}{r(p-k+1)}\|f\|_{L^s(L^r)}M^\frac{(k-1)r-p}{r(p-k+1)}\, t^\frac{s-1}{s} \Bigr]^\frac{r(p-k+1)}{(k-1)(r-1)}
\ee
holding for any $p$ from interval \rf{p_constr2}.

Finally, choosing
\bes
l=\tfrac{2p}{p-k}\in(2,\infty),  \text{ i.e.  so that}\quad p=\tfrac{kl}{l-2},
\ees
in estimate \rf{diss_est} together with \rf{p_bound2} implies
\bes
g(t)\int\limits_{\Omega}|\nabla w|^2\,dx\le\int\limits_{\Omega} v^k|\Delta w|^2\,dx
\ees
with
\be
\lb{g}
g(t)= C_{S,\,\frac{2p}{p-k}}^{-2} \Bigl[\Bigl(\int\limits_{\Omega}\tfrac{dx}{v_0^p}\Bigr)^\frac{(k-1)(r-1)}{r(p-k+1)}
+\tfrac{p(k-1)(r-1)}{r(p-k+1)}\|f\|_{L^s(L^r)}M^\frac{(k-1)r-p}{r(p-k+1)}t^\frac{s-1}{s}\Bigr]^\frac{-kr(p-k+1)}{p(k-1)(r-1)}_.
\ee
Combining that with \rf{w_es2} (while noting that $f_\infty=0$) gives
\bes
\tfrac{d}{dt}\int\limits_{\Omega}|\nabla w|^2\,dx+g(t)\int\limits_{\Omega}|\nabla w|^2\,dx\le
\tfrac{\bar{C}}{\eps^{k/2}}\|f\|_{2+\eps}^2.
\ees
Applying comparison principle to the last inequality one deduces that
\begin{multline*}
 \int\limits_{\Omega}|\nabla w|^2\,dx\le \exp\Bigl\{-\int_0^tg(\tau)\,d\tau\Bigr\}\Big(\int\limits_{\Omega}|\nabla w_0|^2\,dx
+   \\
\tfrac{\bar{C}}{\eps^{k/2}}\int_0^t\exp\Bigl\{\int_0^\tau g(\tilde{\tau})\,d\tilde{\tau}\Bigr\}\|f(\cdot,\tau)\|_{2+\eps}^2\,d\tau\Big).
\end{multline*}
The last estimate would imply asymptotic decay
\be
\lb{as_decay3}
\int\limits_{\Omega}|\nabla w|^2\,dx\go 0\quad\text{as}\quad t\go+\infty
\ee
provided
\be
\lb{g_decay_0}
\int_0^tg(\tau)\,d\tau\go\infty\quad\text{as}\quad t\go\infty,
\ee
as well as,
\bes
\mathop {\lim} \limits_{t\go\infty}\exp\Bigl\{-\int_0^tg(\tau)\,d\tau\Bigr\}\int_0^t\exp\Bigl\{\int_0^\tau g(\tilde{\tau})\,d\tilde{\tau}\Bigr\}\|f(\cdot,\tau)\|_{2+\eps}^2\,d\tau=0
\ees
hold together.
Using L'Hospital's rule and \rf{g}, one observes that the last limit holds if and only if
\be
\lb{f_constr_0}
\mathop {\lim} \limits_{t\go\infty} \|f(\cdot,t)\|_{2+\eps}^2t^\frac{(s-1)kr(p-k+1)}{sp(k-1)(r-1)}=0.
\ee
Similarly, \rf{g_decay_0} holds provided
\bes
\tfrac{(s-1)kr(p-k+1)}{sp(k-1)(r-1)}\le 1.
\ees
The last condition together with the constraint \rf{p_constr2} on $p$  implies necessarily
that
\bes
s\in[1,\,s_{cr}(k,r))\quad\text{with}\quad s_{cr}(k,r)=\tfrac{r}{r-(k-1)(r-1)}.
\ees
Note that $s_{cr}>1$ for $k<k_{cr}(r)$ and $r>1$.

Exploring additionally \rf{f_constr_0} we conclude that for any
\be
\lb{f_cd}
f\in L^s(0,+\infty;L^r(\Omega))
\ee
with $s\in[1,\,s_{cr}(k,r))$ and $r>2$ asymptotic decay \rf{as_decay3} holds.
Substituting $w=v(x,t)-v_\infty(x)$ then yields
\bes
\|\nabla(v(\cdot,t)-v_\infty(.))\|_2\go 0 \text{ and } \|v(\cdot,t)-v_\infty(.) \|_{H^1(\Omega)}\go 0 \text{ as }  t\go+\infty.
\ees

%%%%%%%%%%%%%%%%%%%%%%%%%%%%%%%%%%%%%%%%%%%%%%%%%%%%%%%%%%%%%%%%%
\section{Three-dimensional case}
%%%%%%%%%%%%%%%%%%%%%%%%%%%%%%%%%%%%%%%%%%%%%%%%%%%%%%%%%%%%%%%%%
%%%%%%%%%%%%%%%%%%%%%%%%%%%%%%%%%%%%%%%%%%%%%%%%%%%%%%%%%%%%%%%%%
\subsection{Time-homogeneous case}
%%%%%%%%%%%%%%%%%%%%%%%%%%%%%%%%%%%%%%%%%%%%%%%%%%%%%%%%%%%%%%%%%
In this section, we consider independent of time source function \rf{f},
i.e. $f(x,t)=f_\infty(x)=f(x)$. Given $f\in L^r(\Omega)$ with $r\in(3/2,\infty]$ there exists 
a critical value of exponent 
\be
\lb{k_cr1}
k_{cr}(r)=\tfrac{5r-3}{2r-3}\ge2.5 
\ee
such that for $k>k_{cr}(r)$ we show an exponential asymptotic decay estimate \rf{H1_decay}
for any given positive $|\Omega|,M,\|f\|_r$ and initial condition $v_0(x)\ge0$ provided \rf{ic_as3} holds.
Within the case $k>k_{cr}(r)$ we tactically provide different proofs for two sub-cases $k\in[\frac{4r-2}{r-2},+\infty)$ and $k\in(k_{cr}(r),\,\frac{4r-2}{r-2})$. Showing an unconditional asymptotic decay in the case $k\le k_{cr}(r)$ remains then an open problem.

%%%%%%%%%%%%%%%%%%%%%%%%%%%%%%%%%%%%%%%%%%%%%%%%%%%%%%%%%%%%%%%%%
\subsubsection{Case $k\ge\frac{4r-2}{r-2}>k_{cr}(r)$}
%%%%%%%%%%%%%%%%%%%%%%%%%%%%%%%%%%%%%%%%%%%%%%%%%%%%%%%%%%%%%%%%%
The proof in this case coincides with the one given in subsection 3.1.1 for $N=2$.
We note only that the constraint on power $l$ in estimates \rf{diss_est} and \rf{l} here changes to
\bes
l=\tfrac{2(k-1)(2r-1)}{(k-2)r-k+1}\in(2,6].
\ees
which is satisfied for $k\ge\frac{4r-2}{r-2}$.
Proceeding as in subsection 3.1.1 one obtains decay estimate \rf{H1_decay} with
\bes
\bar{C}=C_{S,\frac{2(k-1)(2r-1)}{(k-2)r-k+1}}^{-2}C^{-\frac{kr}{(k-1)(2r-1)}},
\ees
where $C$ is given by \rf{p_bound}.
Note that one has
\bes
\tfrac{kr}{(k-1)(2r-1)}\go\tfrac{2}{3} \text{ as }  k \go \bigl( \tfrac{4r-2}{r-2}\bigr)^+
\ees
and, therefore, in this limit, asymptotic decay rate $\bar{C}$ tends to a finite number depending on $r$ only.

%%%%%%%%%%%%%%%%%%%%%%%%%%%%%%%%%%%%%%%%%%%%%%%%%%%%%%%%%%%%%%%%%
\subsubsection{Case $\frac{4r-2}{r-2}>k>k_{cr}(r)$}
%%%%%%%%%%%%%%%%%%%%%%%%%%%%%%%%%%%%%%%%%%%%%%%%%%%%%%%%%%%%%%%%%
The proof in this case proceeds similar to the one presented in subsection 2.1.2.
We first choose $p=\frac{3k}{2}$ in estimate \rf{RHS1}, where as before
\bes
k'=\tfrac{(k-1)r}{(r-1)(k-1)-r\beta}
\ees
and $\beta$ is fixed by condition $(p+1-k-\beta)k'=p=3k/2$. This implies
\be
\lb{beta2}
\beta=\tfrac{(2r(k-1)-3k)(k-1)}{r(k+2)}\in(0,k-1)\quad\text{and}\quad k'=\tfrac{(k+2)r}{r(4-k)+2(k-1)}>3,
\ee
where the low bound $k'>3$ follows from the condition $k>k_{cr}(r)$.
Using \rf{RHS1}, \rf{beta2} and Sobolev-Poincare embedding \rf{MTest} taken with $q=6$, and by denoting
\bes
y(t)=\int\limits_{\Omega}\tfrac{dx}{v^\frac{3k}{2}},
\ees
one obtains an analogue of inequality \rf{Entrops3}:
\be
\lb{Entrops5}
y'(t)+\tfrac{6(k+2)}{C_{S,6}^2 \, k}y^\frac{1}{3}(t) \le \tfrac{12(k+2)}{C_{S,6}^2 \, k}\bigl(\tfrac{M}{|\Omega|}\bigr)^\frac{k}{2(k-1)}|\Omega|^\frac{1}{3}
+\tfrac{3k}{2}\|f\|_rM^\frac{k(2r-3)-2r}{r(k+2)}y^\frac{1}{k'}(t).
\ee
Next, using \rf{beta2} and Young's inequality
\bes
ab\le\tfrac{a^s}{s}+\tfrac{b^{s'}}{s'}\quad\text{with}\quad s=\tfrac{k'}{3}\quad\text{and}\quad s'=\tfrac{k'}{k'-3},
\ees
one estimates
\bes
\tfrac{3k}{2}\|f\|_rM^\frac{k(2r-3)-2r}{r(k+2)} y^\frac{1}{k'}(t) \le \tfrac{18(k+2)}{k'\,C_{S,6}^2\,k}y^\frac{1}{3}(t)
+\tfrac{k'-3}{k'}\Bigl(\tfrac{3k}{2}\bigl[\tfrac{C_{S,6}^2k}{6(k+2)}\bigr]^\frac{3}{k'}\|f\|_rM^\frac{k(2r-3)-2r}{r(k+2)}\Bigr)^\frac{k'}{k'-3}_.
\ees
This together with \rf{Entrops5} imply
\begin{multline} \lb{Entrops6}
y'(t)+\tfrac{2(k+2)(k'-3)}{C_{S,6}^2 \, kk'} y^\frac{1}{3}(t)\le   \\
\tfrac{12(k+2)}{C_{S,6}^2k} \bigl(\tfrac{M}{|\Omega|}\bigr)^\frac{k}{2(k-1)}|\Omega|^\frac{1}{3}
+\tfrac{k'-3}{k'}\Bigl(\tfrac{3k}{2}\bigl[\tfrac{C_{S,6}^2k}{6(k+2)}\bigr]^\frac{3}{k'}\|f\|_rM^\frac{k(2r-3)-2r}{r(k+2)}\Bigr)^\frac{k'}{k'-3}.   
\end{multline}
Observe that the right-hand side of \rf{Entrops6} is constant depending only on values of $k,\,r,\,|\Omega|,\,M,\,\|f\|_r$.
Therefore, from \rf{Entrops6} one concludes
\be
\lb{p_bound3}
y(t)=\int\limits_{\Omega}\tfrac{dx}{v^\frac{3k}{2}}\le C=\max\{y(0),C_1\}\quad\forall\, t>0
\ee
with constant
\bes
C_1=\Bigl[\tfrac{6k'}{k'-3}\bigl(\tfrac{M}{|\Omega|}\bigr)^\frac{k}{2(k-1)}|\Omega|^\frac{1}{3}
+\tfrac{C_{S,6}^2k}{2(k+2)}\Bigl(\tfrac{3k}{2}\bigl[\tfrac{C_{S,6}^2k}{6(k+2)}\bigr]^\frac{3}{k'}\|f\|_rM^\frac{k(2r-3)-2r}{r(k+2)}\Bigr)^\frac{3}{k'-3}\Bigr]^3 .
\ees
Finally, choosing $l=6$ in estimate \rf{diss_est} together with \rf{p_bound3} imply \rf{Msc_es} with
\be
\lb{C}
\bar{C}=C_{S,6}^{-2}C^{-\frac{2}{3}},
\ee
and, consequently, decay estimate \rf{H1_decay}.

We note that the derived exponential asymptotic rate $\bar{C}$ does degenerate in the limit $k\go k^+(r)$
as
\bes
k'\go3^+\quad\text{and, consequently,}\quad C\go+\infty \text{ as }  k\go k^+(r).
\ees
\begin{remark}%[Remark]
In the critical case $k=k_{cr}(r)$ one can still show an asymptotic decay result, provided values of $r,\,M$ and $\|f\|_r$
satisfy condition \rf{f_constr} below. In this case, we set $\beta=1,\,k'=3$ in \rf{RHS1} and \rf{beta2}. Then \rf{Entrops5} becomes
\be
\lb{Entrops7}
y'(t)\le\tfrac{12(k_{cr}+2)}{C_{S,6}^2k_{cr}} \bigl(\tfrac{M}{|\Omega|}\bigr)^\frac{k_{cr}}{2(k_{cr}-1)}|\Omega|^\frac{1}{3}
-\Bigl(\tfrac{6(k_{cr}+2)}{C_{S,6}^2k_{cr}}-\tfrac{3k_{cr}}{2}\|f\|_r M^\frac{1}{k_{cr}-1}\Bigr)y^\frac{1}{3}(t) .
\ee
If
\be
\lb{f_constr}
\|f\|_rM^\frac{2r-3}{3r}<\tfrac{4(k_{cr}+2)}{C_{S,6}^2k_{cr}^2}=\tfrac{36}{C_{S,6}^2}\tfrac{(r-1)(2r-3)}{(5r-3)^2},
\ee
one obtains from \rf{Entrops7} that
\bes
y(t)=\int\limits_{\Omega}\tfrac{dx}{v^\frac{3k_{cr}}{2}}\le C=\max\{y(0),C_2\}\quad\forall\, t>0
\ees
with constant
\bes
C_2=\Bigl[\tfrac{8(k_{cr}+2)k_{cr}\bigl(\tfrac{M}{|\Omega|}\bigr)^\frac{k_{cr}}{2(k_{cr}-1)}|\Omega|^\frac{1}{3}}
{4(k_{cr}+2)-C_{S,6}^2k_{cr}^2\|f\|_r M^\frac{1}{k_{cr}-1}}\Bigr]^3 .
\ees
Consequently, choosing $l=6$ in estimate \rf{diss_est} yields decay\rf{H1_decay} with $\bar{C}$ given by \rf{C}.
\end{remark}

%%%%%%%%%%%%%%%%%%%%%%%%%%%%%%%%%%%%%%%%%%%%%%%%%%%%%%%%%%%%%%%%%
\subsection{Time-inhomogeneous case}
%%%%%%%%%%%%%%%%%%%%%%%%%%%%%%%%%%%%%%%%%%%%%%%%%%%%%%%%%%%%%%%%%
In this section, we consider time-dependent source function $f(x,\,t)$ in \rf{f} 
and separate cases $k\in(k_{cr}(r),\,+\infty)$ and $k\in[\frac{r}{r-3/2},\,k_{cr}(r)]$ with $k_{cr}(r)$ given in \rf{k_cr1}. It turns out that in the first case one can take $s=\infty$ in \rf{f} (similar to to the time-homogeneous case $f(x,t)=f(x)$), while in the case $k\in[\frac{r}{r-3/2},\,k_{cr}(r)]$ for the given space integrability exponent $r\in(3/2,\,+\infty]$ we can identify a maximal possible $s_{cr}=s_{cr}(k,r)\in[1,\infty)$, such that asymptotic decay to the steady state $v_\infty(x)$ is shown for any time integrability exponent $s\in[1,\,s_{cr}(k,r)]$ in \rf{f}. For $k<\frac{r}{r-3/2}$ possibility for an asymptotic decay to  $v_\infty(x)$ remains then an open problem. 

Our starting point is given by estimate \rf{w_es}. Due to lack of continuous embedding $H^1(\Omega)\hookrightarrow L^p(\Omega)$ for $p>6$ and $N=3$, we estimate the right-hand side of \rf{w_es} in a different way then in subsection 3.2. For $f$ in \rf{f} we introduce the Bogovskii's vector function $\mathbf{B}f(x,\,t)\in\mathbb{R}^3$ having properties
\be
\lb{B_def}
\mathrm{div}\,\mathbf{B}f(x,\,t)=f(x,t)\quad\text{for all}\quad x\in\Omega\quad\text{and}\quad\mathbf{B}f(x,t)=\mathbf{0}\quad\text{on}\quad\partial\Omega,
\ee
for all $t>0$, and such that
\be
\lb{bog_bound}
\|\mathbf{B}f(\cdot,\,t)\|_{W_0^{1,\,r}}\le C_r\|f(\cdot,\,t)\|_r\quad\text{for all}\quad t>0
\ee
with constant $C_r>0$ depending only on $r$. Existence of $\mathbf{B}f$ having these properties was shown in~\cite{Bo79}.
We note, in particular, due to continuous embedding $W^{1,r}(\Omega)\hookrightarrow L^\frac{Nr}{N-r}(\Omega)$ for $r\le N$ that \rf{f_conv} implies for $N\ge2$
\be
\lb{Bog_conv}
\|\mathbf{B}f(\cdot,\,t)-\mathbf{B}f_\infty(.)\|_2\go 0\quad\text{as}\quad t\go+\infty,
\ee
while \rf{f_tau} implies
\be
\lb{Bf_tau}
\|\mathbf{B}(f_t)(.,t) \|_2\go 0\quad\text{as}\quad t\go+\infty.
\ee
Next, let us rewrite formula \rf{w_es} in the form
\bes
\tfrac{1}{2}\tfrac{d}{dt}\int\limits_{\Omega}|\nabla w|^2\,dx+\tfrac{d}{dt}\int\limits_{\Omega}v(f-f_\infty)\,dx+
\int\limits_{\Omega}v^k(\Delta w-f+f_\infty)^2\,dx=\int\limits_{\Omega}vf_t\,dx.
\ees
Using integration by parts and \rf{B_def} the last identity transforms as
\begin{multline*}
\tfrac{1}{2}\tfrac{d}{dt}\int\limits_{\Omega}|\nabla w-\mathbf{B}(f-f_\infty)|^2\,dx-
\tfrac{1}{2}\tfrac{d}{dt}\int\limits_{\Omega}|\mathbf{B}(f-f_\infty)|^2\,dx +    \\
\int\limits_{\Omega}v^k(\nabla\cdot[\nabla w-\mathbf{B}(f-f_\infty)])^2\,dx   =\int\limits_{\Omega}wf_t\,dx
=-\int\limits_{\Omega}\nabla w\cdot\mathbf{B}(f_t)\,dx= \\
-\int\limits_{\Omega}(\nabla w-\mathbf{B}(f-f_\infty))\cdot\mathbf{B}(f_t)\,dx-\int\limits_{\Omega}\mathbf{B}(f-f_\infty)\cdot\mathbf{B}(f_t)\,dx.
\end{multline*}
This yields
\bes
\tfrac{1}{2}\tfrac{d}{dt}\int\limits_{\Omega}z^2\,dx+\int\limits_{\Omega}v^k(\nabla\cdot z)^2\,dx=-\int\limits_{\Omega}z\cdot\mathbf{B}(f_t)\,dx
\ees
with $z:=\nabla w-\mathbf{B}(f-f_\infty) $. Additionally, arguing as in the proof of estimate \rf{diss_est}, one shows
\bes
\int\limits_{\Omega}z^2\,dx\le C_{S,l}^2\Bigl(\int\limits_{\Omega} v^k|\nabla z|^2\,dx\Bigr)\Bigl(\|\tfrac{1}{v}\|_\frac{kl}{l-2}\Bigr)^k.
\ees
Finally, combining the last two inequalities and using Cauchy-Schwarz inequality yields
\be
\lb{zest}
\tfrac{d}{dt}\int\limits_{\Omega}z^2\,dx+\Bigl(C_{S,l}^{2k}\|\tfrac{1}{v}\|_\frac{kl}{l-2}\Bigr)^{-k}\int\limits_{\Omega}z^2\,dx\le\Bigl(C_{S,l}^{2k}\|\tfrac{1}{v}\|_\frac{kl}{l-2}\Bigr)^k\|\mathbf{B}(f_t)\|_2^2.
\ee
In subsections 4.2.1 and 4.2.2 below, we derive uniform in time bounds on $\|\frac{1}{v}\|_\frac{kl}{l-2}$ for certain $l>2$. These bounds combined with \rf{zest} will then imply the asymptotic decay \rf{non_stat_conv}.

%%%%%%%%%%%%%%%%%%%%%%%%%%%%%%%%%%%%%%%%%%%%%%%%%%%%%%%%%%%%%%%%%
\subsubsection{Case $k>k_{cr}(r)$}
%%%%%%%%%%%%%%%%%%%%%%%%%%%%%%%%%%%%%%%%%%%%%%%%%%%%%%%%%%%%%%%%%

Here, w.l.o.g., we assume that $s=\infty$ and $r\in(3/2,\,+\infty]$ in \rf{f}. In the first part of the proof,
 we follow exactly the lines of subsection 4.1.1 (in the case $k\in[\frac{4r-2}{r-2},\,+\infty)$) or subsection 4.1.2 (in the case $k\in(\frac{4r-2}{r-2},\,k_{cr}(r))$) establishing uniform in time bounds \rf{p_bound} or \rf{p_bound3}, receptively, i.e.
\bes
\int\limits_{\Omega}\tfrac{dx}{v^p}\le C\quad\forall\, t>0 \text{ and fixed }  p\in[\tfrac{3k}{2},\,+\infty),\ C>0.
\ees
We note just that by deriving \rf{p_bound} or \rf{p_bound3} $\|f\|_r$ should be now replaced by $\|f\|_{L^\infty(0,+\infty;L^r(\Omega))}$.

Next, the latter bound combined with \rf{zest} taken with $l=\frac{2p}{p-k}$ imply
\bes
\lb{z_est}
\tfrac{d}{dt}\int\limits_{\Omega}z^2\,dx+C_{S,\frac{2p}{p-k}}^{-2}C^{-\frac{k}{p}}\int\limits_{\Omega}z^2\,dx\le C_{S,\frac{2p}{p-k}}^2C^\frac{k}{p}\|\mathbf{B}(f_t)\|_2^2.
\ees
Applying comparison principle to the last inequality one deduces that
\beas
\int\limits_{\Omega}|\nabla z|^2\,dx&\le& \exp\Bigl\{-C_{S,\frac{2p}{p-k}}^{-2}C^{-\frac{k}{p}}\,t\Bigr\}\Big(\int\limits_{\Omega}|\nabla z_0|^2\,dx+\nonumber\\
&+&C_{S,\frac{2p}{p-k}}^2C^\frac{k}{p}\int_0^t\exp\Bigl\{C_{S,\frac{2p}{p-k}}^{-2}C^{-\frac{k}{p}}\,\tau\Bigr\}\|\mathbf{B}(f_t)\|_2^2\,d\tau\Big).
\eeas
The last estimate would imply asymptotic decay
\be
\lb{as_decayz}
\int\limits_{\Omega}|\nabla z|^2\,dx\go 0\quad\text{as}\quad t\go+\infty 
\ee
provided
\bes
\mathop {\lim} \limits_{t\go\infty}\exp\Bigl\{-C_{S,\frac{2p}{p-k}}^{-2}C^{-\frac{k}{p}}\,t\Bigr\}\int_0^t\exp\Bigl\{C_{S,\frac{2p}{p-k}}^{-2}C^{-\frac{k}{p}}\,\tau\Bigr\}\|\mathbf{B}(f_t)\|_2^2\,d\tau=0.
\ees
Using L'Hospital's rule one observes that this holds because of \rf{Bf_tau}. Hence, \rf{as_decayz} holds and, after substituting $z=\nabla w-\mathbf{B}(f-f_\infty)$ and $w=v(x,t)-v_\infty(x)$, one obtains
\bes
\|\nabla(v(\cdot,t)-v_\infty(.) )\|_2\go 0 \text{ and } \|v(\cdot,t)-v_\infty(.) \|_{H^1(\Omega)}\go 0 \text{ as }  t\go+\infty,
\ees
where we also use \rf{Bog_conv}.

%%%%%%%%%%%%%%%%%%%%%%%%%%%%%%%%%%%%%%%%%%%%%%%%%%%%%%%%%%%%%%%%%
\subsubsection{Case $\frac{r}{r-3/2}\le k\le k_{cr}(r)$}
%%%%%%%%%%%%%%%%%%%%%%%%%%%%%%%%%%%%%%%%%%%%%%%%%%%%%%%%%%%%%%%%%

Here, we assume $f\in L^s(0,+\infty;L^r(\Omega))$ with $s\in[1,\infty)$ and  $r\in(3/2,\,+\infty)$ and argue similarly as in subsection 3.2.2. First, we set $p=\frac{3k}{2}$ and consider estimate \rf{RHS2} with $k'$ and $\beta$ given by \rf{k'} and \rf{beta1}, respectively. Note that both conditions $\beta\in[0,\,k-1)$ and $p+1-k-\beta>0$ hold for our range of $k$. Then from  \rf{Entrops} and \rf{RHS2}--\rf{beta1} considered with  $p=3k/2$ one gets 
\bes
\tfrac{d}{dt}\int\limits_{\Omega}\tfrac{dx}{v^\frac{3k}{2}}\le \tfrac{3k}{2}\|f\|_r M^\frac{2(k-1)r-3k}{r(k+2)}\Bigl(\int\limits_{\Omega}\tfrac{dx}{v^p}\Bigr)^\frac{3kr-2(k-1)(2r-1)}{r(k+2)} .
\ees
Next, applying Gr\"{o}nwall's inequality to this estimate results in the bound
\be
\lb{p_bound4}
\int\limits_{\Omega}\tfrac{dx}{v^\frac{3k}{2}}\le \Bigl[\Bigl(\int\limits_{\Omega}\tfrac{dx}{v_0^\frac{3k}{2}}\Bigr)^\frac{2(k-1)(r-1)}{r(k+2)}
+\tfrac{3k(k-1)(r-1)}{r(k+2)}\|f\|_{L^s(L^r)}M^\frac{2(k-1)r-3k}{r(k+2)}\, t^\frac{s-1}{s} \Bigr]^ \frac{r(k+2)}{2(k-1)(r-1)}.
\ee
Choosing
\bes
l=6, \text{ i.e.  so that }  p=\tfrac{kl}{l-2},
\ees
in estimate \rf{zest} together with \rf{p_bound4} imply
\bes
\tfrac{d}{dt}\int\limits_{\Omega}z^2\,dx+g(t)\int\limits_{\Omega}z^2\,dx\le\tfrac{\|\mathbf{B}(f_t)\|_2^2}{g(t)}.
\ees
with
\be
\lb{g1}
g(t)= C_{S,\,6}^{-2} \Bigl[\Bigl(\int\limits_{\Omega}\tfrac{dx}{v_0^p}\Bigr)^\frac{2(k-1)(r-1)}{r(k+2)}
\!\!\!\! +\tfrac{3k(k-1)(r-1)}{r(k+2)}\|f\|_{L^s(L^r)}M^\frac{2(k-1)r-3k}{r(k+2)}t^\frac{s-1}{s}\Bigr]^ \frac{-r(k+2)}{3(k-1)(r-1)} 
\! \!.
\ee
Applying comparison principle to the last inequality one deduces that
\begin{multline*}
 \int\limits_{\Omega}|\nabla z|^2\,dx\le \exp\Bigl\{-\int_0^tg(\tau)\,d\tau\Bigr\} \times  \\
\Bigl(\int\limits_{\Omega}|\nabla z_0|^2\,dx
+\int_0^t\exp\Bigl\{\int_0^\tau g(\tilde{\tau})\,d\tilde{\tau}\Bigr\}\tfrac{\|\mathbf{B}(f_\tau)\|_2^2}{g(\tau)}\,d\tau\Bigr).   
\end{multline*}
The last estimate would imply asymptotic decay
\be
\lb{as_decay5}
\int\limits_{\Omega}|\nabla z|^2\,dx\go 0\quad\text{as}\quad t\go+\infty 
\ee
provided  
\be
\lb{g_decay}
\int_0^tg(\tau)\,d\tau\go\infty\quad\text{as}\quad t\go\infty,
\ee
as well as,
\bes
\mathop {\lim} \limits_{t\go\infty}\exp\Bigl\{-\int_0^tg(\tau)\,d\tau\Bigr\}\int_0^t\exp\Bigl\{\int_0^\tau g(\tilde{\tau})\,d\tilde{\tau}\Bigr\}\tfrac{\|\mathbf{B}(f_t)\|_2^2}{g(t)}\,d\tau=0
\ees
hold together.
Using L'Hospital's rule and \rf{g1}  one observes that the last limit holds if and only if
\be
\lb{f_constr1}
\mathop {\lim} \limits_{t\go\infty} \|\mathbf{B}(f_t)\|_2\cdot t^\frac{(s-1)r(k+2)}{6s(k-1)(r-1)}=0.
\ee
Similarly, \rf{g_decay} holds provided
\bes
\tfrac{(s-1)r(k+2)}{3s(k-1)(r-1)}\le 1.
\ees
The last condition implies, necessarily,
that
\bes
s\in[1,\,s_{cr}(k,r)]\quad\text{with}\quad s_{cr}(k,r)=\tfrac{r(k+2)}{r(5-2k)+3(k-1)}.
\ees
Note that $s_{cr}(k,r)>1$ for $k<k_{cr}(r)$ and $r>1$.

Exploring additionally \rf{f_constr1} and using \rf{bog_bound}, \rf{f_tau} we conclude that for any
\bes
f\in L^s(0,+\infty;L^r(\Omega))
\ees
with $s\in[1,\,s_{cr}(k,r)]$ decay \rf{as_decay5} holds.
Substituting there\\ $z=\nabla w-\mathbf{B}(f-f_\infty)$ and $w=v(x,t)-v_\infty(x)$, one obtains
\bes
\|\nabla(v(\cdot,t)-v_\infty(.))\|_2\go 0 \text{ and } \|v(\cdot,t)-v_\infty(.)\|_{H^1(\Omega)}\go 0 \text{ as }  t\go+\infty,
\ees
where we also use \rf{Bog_conv}.

%%%%%%%%%%%%%%%%%%%%%%%%%%%%%%%%%%%%%%%%%%%%%%%%%%%%%%%%%%%%%%%%%
\section{Discussion}
%%%%%%%%%%%%%%%%%%%%%%%%%%%%%%%%%%%%%%%%%%%%%%%%%%%%%%%%%%%%%%%%%
In this study, we identified critical mobilities $k=k_{cr}(r)$ \rf{k_cr} and \rf{k_cr1} for problem \rf{ME1}--\rf{ME3} considered with source term \rf{f} in dimensions $N=2$ and $N=3$. In dimension $N=2$, Theorems 1.1--1.2 imply asymptotic decay in $H^1$-norm of dynamical solutions to the unique positive state $u_\infty$ for $k\ge k_{cr}(r)$ and provide with explicit decay rates. In dimension $N=3$, the same is shown for $k>k_{cr}(r)$ and possibility of asymptotic decay for  $k=k_{cr}(r)$ remains then an open question. Interestingly, condition for existence of a unique positive steady state $v_\infty$ has also form $k\ge k_{cr}(r)$ for $r\in(N/2,\,+\infty]$ (see Appendix A.2). This suggests that for $k<k_{cr}(r)$ possibility of decay to a nonhomogeneous steady state may be conditional to certain relations between given $M$,$|\Omega|$ and $\|f\|_r$.
Besides that, possibility of decay to constant profile $v_\infty$ corresponding to $f_\infty=0$ (see \rf{const_ss} in Appendix A.2) for $k\in(1,\,\frac{r}{r-1}]$ for $N=2$ and  $k\in(1,\,\frac{r}{r-3/2})$ for $N=3$ is unknown.

Other open questions related to the assumptions of Theorems 1.1-1.2 include question about positivity of solutions for all $t>0$ for positive initial profiles $v_0>0$ in case $r\in(N/2,\,+\infty)$ as well as existence and regularity properties of solutions starting from nonnegative initial data $v_0\ge 0$. Possible existence and decay of sign-changing solutions to \rf{ME1}--\rf{ME3}  considered with \rf{f} in certain rages of $k>1$ and $p=1$ remained also out of the scope of this study.

Finally, investigation of the limiting cases $k\go 1^+$ or $k\go+\infty$ in \rf{ME1}--\rf{ME3} corresponding to $m\go+\infty$ and $m\go 1$ in the PME problem \rf{PME1}--\rf{PME3} might be of interest. In both of the limits
many of our estimates degenerate.

\begin{appendix}
%%%%%%%%%%%%%%%%%%%%%%%%%%%%%%%%%%%%%%%%%%%%%%%%%%%%%%%%%%%%%%%%%
\section{Positivity of solutions}
%%%%%%%%%%%%%%%%%%%%%%%%%%%%%%%%%%%%%%%%%%%%%%%%%%%%%%%%%%%%%%%%%
In this section, we discuss positivity properties of solutions  $v(x,t)$ to problem \rf{ME1}--\rf{ME3} and
of the stationary one $v_\infty(x)$ to \rf{SME1}--\rf{SME3}.

%%%%%%%%%%%%%%%%%%%%%%%%%%%%%%%%%%%%%%%%%%%%%%%%%%%%%%%%%%%%%%%%%
\subsection{Positivity of dynamic solutions}
%%%%%%%%%%%%%%%%%%%%%%%%%%%%%%%%%%%%%%%%%%%%%%%%%%%%%%%%%%%%%%%%%
Using a priori estimate \rf{Entrops}, we can show that solutions $v(x,t)$ to \rf{ME1}--\rf{ME3} stay positive for all $t>0$ and $x\in\bar{\Omega}$  provided initially $v_0(x)>0$. This result is shown here for all parameter values $k\in(1,\,+\infty)$ and in any dimension $N\ge 1$, but only in case $r=+\infty$ considered in \rf{f}. Positivity of solutions in cases  $r\in(N/2,\,+\infty)$ then remains an open problem.

The proof in case $r=+\infty$ proceeds as follows. Let us estimate the term at the right-hand side of \rf{Entrops} as
\bes
\int\limits_{\Omega}\tfrac{f\,dx}{v^{p+1-k}}\le\|f\|_\infty\Bigl(\int\limits_{\Omega}\tfrac{dx}{v^p}\Bigr)^\frac{p+1-k}{p}|\Omega|^\frac{k-1}{p},
\ees
and denote
\bes
y(t)=\int\limits_{\Omega}\tfrac{dx}{v^p}.
\ees
Then \rf{Entrops} implies
\bes
y'(t)\le p\|f\|_\infty|\Omega|^\frac{k-1}{p}y^\frac{p+1-k}{p}(t).
\ees
Applying Gr\"{o}nwall's inequality to this estimate results in the following bound
\bes
y(t)\le\Bigl(y(0)^\frac{k-1}{p}+(k-1)\|f\|_{L^s(0,+\infty;L^\infty(\Omega))}|\Omega|^\frac{k-1}{p}t^\frac{s}{s-1}\Bigr)^\frac{p}{k-1},
\ees
i.e.
\bes
\|v(\cdot,\,t)^{-1}\|_p\le\Bigl(\|v_0^{-1}\|_p^{k-1}+(k-1)\|f\|_{L^s(0,+\infty;L^\infty(\Omega))}|\Omega|^\frac{k-1}{p}t^\frac{s}{s-1}\Bigr)^\frac{1}{k-1}.
\ees
Letting   $p\go\infty$ in this inequality gives	
\bes
\|v(\cdot,\,t)^{-1}\|_\infty\le\Bigl(\|v_0^{-1}\|_\infty^{k-1}+(k-1)\|f\|_{L^s(0,+\infty;L^\infty(\Omega))}t^\frac{s}{s-1}\Bigr)^\frac{1}{k-1}.
\ees
Hence, $v(x,t)>0$ for all $t>0$ and $x\in\bar{\Omega}$  provided $v_0>0$. We note that the obtained bound from below on
$v(x,t)$ is not uniform in time, but still implies smoothness and uniqueness of solutions $v(\cdot,t)>0$ for all $t>0$.

%%%%%%%%%%%%%%%%%%%%%%%%%%%%%%%%%%%%%%%%%%%%%%%%%%%%%%%%%%%%%%%%%
\subsection{Existence of positive steady states $v_\infty(x)$}
%%%%%%%%%%%%%%%%%%%%%%%%%%%%%%%%%%%%%%%%%%%%%%%%%%%%%%%%%%%%%%%%%
Let us recall the definitions \rf{k_cr} and \rf{k_cr1} for the critical value of exponent $k$ in dimensions $N=2$ and $N=3$, respectively. First, we show that for any parameter $k\in [k_{cr}(r),\,+\infty)$ there exists a unique positive stationary solution $v_\infty$ to problem  \rf{SME1}--\rf{SME3}, provided the limit right-hand side $f_\infty\in L^r(\Omega)$. This follows from conservation  \rf{SME3} combined with Sobolev continuous embedding (see, e.\,g.~\cite{AF03}) 
\be
\lb{embed}
W^{2,r}(\Omega)\hookrightarrow C^\lambda(\bar{\Omega})\quad\text{for}\quad0<\lambda \le 2-\tfrac{N}{r}  \text{ and } r > \tfrac{N}{2}.
\ee
Let $v_\infty$ be a solution to \rf{SME1}--\rf{SME2}. The embedding implies that
\bes
|v_\infty(x)|\le L|x-x_0|^\lambda
\ees
holds for any $x\in\bar{\Omega}$ and some constant $L>0$ whenever $v_\infty(x_0)=0$.
Let us further consider nonnegative solutions $v_\infty$ to \rf{SME1}--\rf{SME2}.
For them \rf{SME3} implies that
\be
\lb{hol}
M=\int\limits_{\Omega}\tfrac{dx}{v_\infty^{k-1}}\ge \frac{1}{L^{k-1}}\int\limits_{\Omega}\tfrac{dx}{|x-x_0|^{\lambda(k-1)}}
\ee
The integral at the right-hand side of \rf{hol} diverges provided $\lambda(k-1)\ge N$, i.e. using \rf{embed} when
\bes
k > \tfrac{N+ 2}{2} \text{ and } r\ge \tfrac{N(k-1)}{2(k-1)-N}
\ees
holds. We note that the latter condition coincides with $k\in [k_{cr}(r),\,+\infty)$ for $N=2,\,3$ and, consequently, condition $f_\infty\in L^r(\Omega)$ implies existence of a unique positive solution  $v_\infty$ to \rf{SME1}--\rf{SME3}.

In case $k\in(1,\,k_{cr}(r))$, one has $f_\infty = 0$ in Theorem 1.2 and there exists unique positive solution $v_\infty$ to \rf{SME1}--\rf{SME3}:
\be
\lb{const_ss}
v_\infty=\bigl(\tfrac{|\Omega|}{M}\bigr)^\frac{1}{k-1}.
\ee
We conclude that assumptions of Theorems 1.1-1.2 imply, in particular, existence of unique positive $v_\infty$.

%%%%%%%%%%%%%%%%%%%%%%%%%%%%%%%%%%%%%%%%%%%%%%%%%%%%%%%%%%%%%%%%%
\section{Boundedness of $\|\nabla v\|_2$}
%%%%%%%%%%%%%%%%%%%%%%%%%%%%%%%%%%%%%%%%%%%%%%%%%%%%%%%%%%%%%%%%%
Here, we prove formula \rf{grad_v} provided \rf{f_tau}--\rf{f_conv} hold. From energy inequality \rf{EE} one obtains
\beas
\tfrac{1}{2}\|\nabla v\|_2^2&\le&E(v_0)-\int\limits_{\Omega}fv\,dx+ \iint \limits_{Q_t} {vf_t\,dx\,d\tau}\nonumber\\
&\le&E(v_0) + \|v\|_2\|f\|_2+\int_0^t\|v\|_2\|f_\tau\|_2\,d\tau.
\eeas
Combining this estimate with Poincare inequality
\bes
\|v\|_2-\Bigr(\tfrac{|\Omega|}{M}\Bigl)^\frac{1}{k-1}|\Omega|^\frac{1}{2}\le\|v-\Bigr(\tfrac{|\Omega|}{M}\Bigl)^\frac{1}{k-1}\|_2\le C_P\|\nabla v\|_2,
\ees
gives
\bea
\lb{Msc1}
\tfrac{1}{2}\|\nabla v\|_2^2&\le&E(v_0)+\|f\|_2\Bigl(C_P\|\nabla v\|_2+\Bigr(\tfrac{|\Omega|}{M}\Bigl)^\frac{1}{k-1}|\Omega|^\frac{1}{2}\Bigr)\nonumber\\
&+&\int_0^t\Bigr(C_P\|\nabla v\|_2+ \bigl(\tfrac{|\Omega|}{M}\bigl)^\frac{1}{k-1}|\Omega|^\frac{1}{2}\Bigr)\|f_\tau\|_2\,d\tau.
\eea
Additionally, one has
\bes
\tfrac{d}{dt}\|f\|_2=\tfrac{<f,\,f_t>}{\|f\|_2}\le\|f_t\|_2,
\ees
i.e.
\bes
\|f\|_2\le\|f_0\|_2+\int_0^t\|f_\tau\|_2\,d\tau.
\ees
Combining the last estimate with \rf{Msc1}, while denoting
\bes
m(t)=\int_0^t\|f_\tau\|_2\,d\tau\quad\text{and}\quad m_\infty=\int_0^\infty\|f_\tau\|_2\,d\tau<\infty,
\ees
yields
\begin{multline}\lb{Msc2} 
(\|\nabla v\|_2-C_P\|f_0\|_2-C_Pm_\infty)^2 \le  2E(v_0)+C_P^2(\|f_0\|_2+m_\infty)^2+ \\
2C_P\int_0^t\|\nabla v\|_2m'(\tau)\,d\tau  + 2(\|f_0\|_2+2m_\infty)\Bigl(\tfrac{|\Omega|}{M}\Bigr)^\frac{1}{k-1}|\Omega|^\frac{1}{2} \\
\le a_0+2\int_0^t\Bigl|\|\nabla v\|_2-C_P\|f_0\|_2-C_Pm_\infty\Bigr|m'(\tau)\,d\tau,
\end{multline}
where we denoted
\begin{multline*}
 a_0=2E(v_0)+C_P^2(\|f_0\|_2+m_\infty)^2+2m_\infty C_P(\|f_0\|_2+   \\
2m_\infty)+2(\|f_0\|_2+2m_\infty)\Bigl(\tfrac{|\Omega|}{M}\Bigl)^\frac{1}{k-1}|\Omega|^\frac{1}{2}.   
\end{multline*}
The last inequality after an application of Bihari-LaSalle lemma to\\  $(\|\nabla v\|_2-\|f_0\|_2-m_\infty)^2$, while
using $m'(t)\ge0$ and $m(0)=0$, implies
\bes
|\|\nabla v\|_2-C_P\|f_0\|_2-C_Pm_\infty|\le\sqrt{2a_0}+m(t)\le\sqrt{2a_0}+m_\infty,
\ees
whence we obtain the sought estimate \rf{grad_v} holding for all $t>0$ with
\bes
C=\sqrt{2a_0}+(C_P+1)m_\infty+C_P\|f_0\|_2.
\ees

\end{appendix}


\begin{thebibliography}{1}
	
	\bibitem{DP97}
	P. Daskalopoulus and M. del Pino.
	\newblock On nonlinear parabolic equations of very fast diffusion.
	\newblock {\em  Arch. Rational Mech. Anal.}, 137: 363-380, 1997.
		
	\bibitem{Va07}
	J.~L. Vazquez.
	\newblock Porous medium equation. Mathematical theory.
	\newblock Oxford mathematical monographs.
	\newblock Claredon Press, Oxford 2007.
	
	\bibitem{Va06}
	J.~L. Vazquez.
	\newblock Smoothing and decay estimates for nonlinear diffusion equations. Equations of porous medium type.
	\newblock Oxford lecture series in Mathematics and Applications 33.
	\newblock Oxford University Press, Oxford 2006.	
		
	\bibitem{GM13}
	G. Grillo and M. Muratori.
	\newblock Sharp short and long time $L^\infty$ bounds for solutions to porous media equations with homogeneous Neumann boundary conditions.
	\newblock {\em  J.Differential Equations}, 254: 2261-2288, 2013.	
	
	\bibitem{SSZ16}
	A. Segatti, G. Schimperna and S. Zelik.
	\newblock On a singular heat equation with dynamic boundary conditions.
	\newblock {\em  Asymptotic Analysis}, 97: 27-59, 2016.
	
	\bibitem{FKT18}
	M.~A. Fontelos, G. Kitavtsev and R.~M. Taranets.
	\newblock Asymptotic decay and non-rupture of viscous sheets.
	\newblock {\em  Z. Angew. Math. Phys.}, 69: 79-, 2018.
		
	\bibitem{KT20}
	G. Kitavtsev and R.~M. Taranets.
	\newblock Long-time behavior of solutions to a singular heat equation with an application to hydrodynamics.
	\newblock {\em  J. Interfaces and Free Boundaries}, 22(2): 157-174, 2020.
	
	\bibitem{Br09}
	D. Bresch.
	\newblock Shallow water equations and related topics.
	\newblock In Handbook of Differential equations. Evolutionary equations volume 5, 1-104
	\newblock Elsevier 2009.
	
    \bibitem{KLN11}
    G. Kitavtsev, P. Laurencot and B. Niethammer.
    \newblock Weak solutions to lubrication equations in presence of strong slippage.
    \newblock {\em  Methods Appl. Anal.}, 18: 183-202, 2011.
        
    \bibitem{PF90}
    O. Penrose and P.~C. Fife.
    \newblock Thermodynamically consistent models of phase-field type for the kinetics of phase transitions.
    \newblock {\em  Phys. D}, 43: 44-62, 1990.
	
	\bibitem{MTSO17}
	M. Mizuguchi, K. Tanaka, K. Sekine and S. Oishi.
	\newblock Estimation of Sobolev embedding constant on a domain dividable into bounded convex domains.
	\newblock {\em  J. Inequal. Appl.}, 299, 2017.
	
	\bibitem{LL01}
	E.~H. Lieb and M. Loss.
	\newblock Analysis.
	\newblock  Second edition. Providence, Rhode Island: American Mathematical Society 2001.
	
	\bibitem{AF03}
	R.~A. Adams and J.~J.~F. Fournier.
	\newblock Sobolev spaces.
	\newblock Pure and applied mathematics series.
	\newblock Second edition. Elsevier 2003.
		
	\bibitem{Bo79}
	M.~E. Bogovskii.
	\newblock Solution of the first boundary value problem for the equation of continuity of an incompressible medium.
	\newblock {\em  Soviet Math. Dokl.}, 20: 1094-1098, 1979.
	
	\bibitem{Bi56}
	I.~A. Bihari.
	\newblock A generalization of a lemma of Bellman and its applications to uniqueness problems of differential equations.
	\newblock {\em  Acta Mathematica Hungarica}, 7: 81-94, 1956.	
	
\end{thebibliography}
\end{document}